\documentclass[11pt]{amsart}
\usepackage{amsmath}
\usepackage[utf8]{inputenc}
\usepackage{amssymb}
\usepackage{amsopn}
\usepackage{epsfig}
\usepackage{amsfonts}
\usepackage{latexsym}
\usepackage{graphicx}
\usepackage{enumerate}
\usepackage{color}
\usepackage{tikz}
\newtheorem{theorem}{Theorem}[section]
\newtheorem{lemma}[theorem]{Lemma}

\newtheorem{corollary}[theorem]{Corollary}

\theoremstyle{definition}

\theoremstyle{remark}
\newtheorem{remark}[theorem]{Remark}

\numberwithin{equation}{section}

\begin{document}
\title{A nonlinear version of the Newhouse thickness theorem}
\author[K. Jiang]{Kan Jiang}
\address[K. Jiang]{Department of Mathematics, Ningbo University,
People's Republic of China}
\email{jiangkan@nbu.edu.cn;kanjiangbunnik@yahoo.com}
\date{\today}
\begin{abstract}
Let $C_1$ and $C_2$ be two Cantor sets with  convex hull $[0,1]$. Newhouse \cite{Newhouse} proved if
$\tau(C_1)\cdot \tau(C_2)\geq 1$, then  the arithmetic sum $C_1+C_2$ is an interval, where $\tau(C_i), 1\leq i\leq 2$ denotes the thickness of $C_i$. In this paper, we  generalize this thickness theorem as follows.
Let $K_i\subset \mathbb{R}, i=1,\cdots, d$, be some Cantor sets  (perfect and nowhere dense)  with  convex hull $[0,1]$. Suppose $f(x_1,\cdots, x_{d-1},z)\in \mathcal{C}^1$ is a continuous function defined on $\mathbb{R}^d$. Denote the continuous image of $f$ by $$f(K_1,\cdots,  K_d)=\{f(x_1, \cdots x_{d-1},z):x_i\in K_i,z\in K_d, 1\leq i\leq d-1\}.$$
If for any $(x_1, \cdots, x_{d-1},z)\in [0,1]^d$, we have
$$(\tau(K_i))^{-1}\leq \left|\dfrac{\partial_{x_i} f}{\partial_z f}\right|\leq \tau(K_d),1\leq i\leq d-1$$ then
$f(K_1,\cdots,  K_d)$ is a closed interval.
We give two applications. Firstly, we partially answer some questions posed by  Takahashi \cite{Takahashi}.  Secondly, we obtain various nonlinear identities, associated with the continued fractions with restricted partial quotients,
 which can represent  real numbers.
\end{abstract}

\maketitle
\section{Introduction}
Let $K_1$ and $K_2$ be  two  Cantor sets with convex hull $[0,1]$. Newhouse \cite{Newhouse} proved that if $\tau(K_1)\cdot \tau(K_2)>1$ (in fact we may replace this condition by $\tau(K_1)\cdot \tau(K_2)\geq 1$),
 then the arithmetic sum $K_1+K_2$ is an interval, where $\tau(K_i)$ denotes the thickness of $K_i, i=1,2.$
 The arithmetic sum of Cantor sets appears naturally in   bifurcation theory.  Palis \cite{Palis} posed the following problem which is currently known as the Palis' conjecture. Whether it is true (at least generically) that the arithmetic  sum of dynamically defined Cantor sets either has measure zero or contains an interior. This conjecture was solved by Moreira and Yoccoz \cite{Yoccoz}. The Newhouse's thickness theorem is a very powerful result which can judge whether the arithmetic sum of two Cantor sets contains interior. Astels \cite[Theorem 2.4]{Astels} generalized the Newhouse's thickness theorem by considering multiple sum of Cantor sets. He made use of this new thickness theorem to prove some identities which can represent real numbers. Astels' thickness theorem implies many interesting results. For instance, we may  prove some Waring type result as follows, see \cite{Tyson, Yu}.
 For each $k\geq 2$, there is a number $n(k)\leq 2^k$ such that for  any  $x\in [0,n(k)]$, we have $$x=\sum_{i=1}^{n(k)}x_i^k,$$ where $x_i$
 is taken from the middle-third Cantor set. 
  The Newhouse and  Astels' thickness theorems are very useful when we consider the sum of two Cantor sets.
  It is natural to consider a nonlinear version of Newhouse's thickness theorem.   Suppose $f(x_1,\cdots, x_{d-1},z)\in \mathcal{C}^1$ is a continuous function defined on $\mathbb{R}^d$. Denote the continuous image of $f$ by $$f(K_1,\cdots,  K_d)=\{f(x_1, \cdots, x_{d-1},z):x_i\in K_i,z\in K_d, 1\leq i\leq d-1\},$$ where $\{K_i\}_{i=1}^{d}$  are  general Cantor sets. To the best of our knowledge, there are very few results about $f(K_1,\cdots,  K_d)$.
Generally, to consider the topological structure of $f(K_1,\cdots,  K_d)$ is a difficult question. As we know very little information about $\{K_i\}_{i=1}^{d}$.  Moreover, the nonlinearity of $f(x_1, \cdots, x_{d-1},z)$ makes the abstract set $f(K_1,\cdots,  K_d)$ obscure.  The main aim of this paper is to give some sufficient conditions on $f(x_1, \cdots, x_{d-1},z)$ such that  $f(K_1,\cdots,  K_d)$ is a closed interval.

We now introduce some related results concerning with the continuous image of $f$ in $\mathbb{R}$. The first one, to the best of our knowledge, is due to
Steinhaus \cite{HS} who  proved in 1917 the following interesting  results:
$$C+C=\{x+y:x,y\in C\}=[0,2], C-C=\{x-y:x,y\in C\}=[-1,1],$$ where $C$ is the middle-third Cantor set.
It  is worth pointing out that Steinhaus also proved that for any two sets with positive Lebesgue measure, their arithmetic sum contains interiors.
In 2019, Athreya, Reznick and Tyson \cite{Tyson} proved that
\[
C\div C=\left\{\dfrac{x}{y}:x,y\in C, y\neq0\right\}=\bigcup_{n=-\infty}^{\infty}\left[ 3^{-n}\dfrac{2}{3},3^{-n}\dfrac
{3}{2}\right] \cup\{0\}.
\]
In \cite{Gu}, Gu, Jiang, Xi and Zhao gave the topological structure of $$C\cdot C=\{xy:x,y\in C\}.$$
They proved that the exact Lebesgue measure of $C\cdot C$ is about $0.80955.$  We give some remarks on the above results. The main idea of \cite{Tyson} is effective for homogeneous self-similar sets. For a general self-similar set or some general Cantor set, we may not utilize their idea directly.
Fraser,  Howroyd and Yu \cite{Fraser} studied the dimensions of sumsets and iterated sumsets, and provided natural conditions which guarantee that a set $F\subset \mathbb{R}$ satisfies $\overline{\dim_{B}}F+F>\overline{\dim_{B}}F$. The reader can find more related references in \cite{Fraser}.

For  higher dimensions, namely $\mathbb{R}^d , d\geq 3$, there are  relatively few results. Banakh, Jab{\l}o\'{n}ska and  Jab{\l}o\'{n}ski \cite{BJJ}  proved under some mild conditions that the arithmetic sum of $d$ many compact connected sets in $\mathbb{R}^d$ has non-empty interior. As a consequence, every compact connected set in $\mathbb{R}^d$ not lying a hyperplane is arithmetically thick.
A compact set $E\subset \mathbb{R}^d$
is said to be arithmetically thick if there exists
a positive integer $n$ so that the $n$-fold arithmetic sum of $E$ has non-empty interior. Recently Feng and Wu \cite{Feng2020} defined the thickness of sets in $ \mathbb{R}^d$, and   proved  the arithmetic thickness
for several classes of fractal sets, including self-similar sets, self-conformal sets in $\mathbb{R}^d$
(with $d \geq 2$) and some self-affine sets. All these elegant  results are concerning with arithmetic sum.
They introduced some new ideas which are very useful to analyze the sets in $ \mathbb{R}^d.$

In this paper,  we consider  similar problems.  However, our main motivation is to  generalize the Newhouse's thickness theorem for some general functions.
Before we introduce the main results of this paper, we give some definitions. First, we give a well-known method that can generate a Cantor set. For simplicity, we let  $I_0=[0,1]$. In the first level, we delete an open interval from $[0,1]$,  denoted by $O$. Then there are two closed intervals left, we denote them by $B_1$ and $B_2$. Therefore, $[0,1]=B_1\cup O\cup B_2.$
Let $E_1=B_1\cup B_2$. In the second level, let  $O_0$ and $O_1$ be open intervals  that are deleted from $B_1$ and $B_2$ respectively, then we clearly have
 $$B_1=B_{11}\cup O_0 \cup B_{12}, B_2=B_{21}\cup O_1 \cup B_{22}.$$ Let $$E_2=B_{11}\cup B_{12}\cup B_{21}\cup B_{22}.$$
Repeating this process, we can generate $E_{n+1}$ from $E_n$ by removing an open interval from each closed interval in the union which consists of $E_n$.  We assume that the deleted open intervals are arranged by the    decreasing lengths, i.e. the  lengths of deleted open intervals are decreasing. If for some levels, the deleted open intervals have the same length, then we can delete these open intervals in any order.
To avoid triviality, we make the following rule.
Let $B_{\omega}$ be a closed interval in some level, then we delete an open interval $O_{\omega}$ from $B_{\omega}$, i.e.
$$B_{\omega}=B_{\omega1}\cup O_{\omega}\cup B_{\omega2}.$$ We assume that the length of $O_{\omega}$ is  positive and strictly smaller than  $B_{\omega}$.
We let $$K=\cap_{n=1}^{\infty}E_n,$$ and call $K$ a Cantor set. The above rule is to ensure the Cantor set  is perfect and nowhere dense.

The next definition is the famous Newhouse's thickness. Given a Cantor set $$K=\cap_{n=1}^{\infty}E_n.$$
Let $B_{\omega}$ be a closed interval in some level. Then by the construction of $K$, we have
$$B_{\omega}=B_{\omega1}\cup O_{\omega}\cup B_{\omega2},$$
where $O_{\omega}$ is an open interval while $B_{\omega1}$ and $B_{\omega2}$ are closed intervals. We call $B_{\omega1}$ and $B_{\omega2}$ bridges of $K$, and $O_{\omega}$ gap of $K$.
Let
$$\tau_{\omega}(B_{\omega})=\min\left\{\dfrac{|B_{\omega1}|}{|O_{\omega}|},\dfrac{|B_{\omega2}|}{|O_{\omega}|}\right\}, $$ where
$|\cdot|$ means length.
We define the thickness of $K$ by
$$\tau(K)=\inf_{B_{\omega}}\tau_{\omega}(B_{\omega}).$$
Here the infimum   takes over all bridges in every level.

Now we  state the  main result of this paper.
\begin{theorem}\label{Main}
Let $\{K_i\}_{i=1}^{d}$  be  Cantor sets with convex hull $[0,1]$. Suppose that $f(x_1,\cdots, x_{d-1},z)\in \mathcal{C}^1$.
If for any $(x_1, \cdots, x_{d-1},z)\in [0,1]^d$, we have
$$(\tau(K_i))^{-1}\leq \left|\dfrac{\partial_{x_i} f}{\partial_z f}\right|\leq \tau(K_d),1\leq i\leq d-1$$ then
$$f(K_1,\cdots,  K_d)=H,$$ where
$$H=\left[\min_{(x_1, \cdots,z)\in K_1\times \cdots  \times K_d}f(x_1,\cdots,z), \max_{(x_1, \cdots,z)\in K_1\times \cdots  \times K_d}f(x_1,\cdots,z)\right],$$ and
$\tau(K_i), i=1,\cdots, d,$ denotes the thickness of $K_i$.
\end{theorem}
\begin{remark}
This result partially generalizes  \cite[Theorem 2.4]{Astels}.
Theorem \ref{Main} can be given an explanation from geometric measure theory. If the convex hull of each $K_i$ is different, then we need to assume each $K_i$ is not contained in any other $K_j$'s gaps, where $j\neq i$. As for this case $f(K_1,\cdots,  K_d)$ may have Lebesgue measure zero.  When we use Theorem \ref{Main}, we need to abide  by this rule.
\end{remark}
\begin{corollary}\label{cor1}
Let $\{K_i\}_{i=1}^{d}$  be  Cantor sets with convex hull $[0,1]$. Suppose that $f(x_1,\cdots x_{d-1},z)\in \mathcal{C}^1$.
If for any $(x_1, \cdots, x_{d-1},z)\in [0,1]^d$, we have
$$(\tau(K_i))^{-1}\leq \left|\dfrac{\partial_{x_i} f}{\partial_z f}\right|\leq \tau(K_d),1\leq i\leq d-1$$
then for any $w\in H$ the hypersurface $f(x_1, \cdots, x_{d-1},z)=w$ intersects with $K_1\times  \cdots\times K_d$.
\end{corollary}
For $d=2$ we have the following result which can be viewed as a nonlinear version of the Newhouse's thickness theorem.
\begin{corollary}\label{cor2}
Let $K_1$ and $K_2$ be two Cantor sets with convex hull $[0,1]$. Suppose $f(x,y)\in \mathcal{C}^{1}.$
 If for any $(x,y)\in [0,1]^2$, we have
$$(\tau(K_1))^{-1}\leq \left|\dfrac{\partial_x f}{\partial_y f}\right|\leq \tau(K_2),$$ then
$$f(K_1,K_2)=\left[\min_{(x,y)\in K_1  \times K_2}f(x,y), \max_{(x,y)\in K_1  \times K_2}f(x,y)\right]=H,$$ where
$\tau(K_i), i=1,2$ denotes the thickness of $K_i$.
In particular, if we take a linear function $$f(x,y)=x+y,  \mbox{ and } \tau(K_i)\geq1, i=1,2,$$ then $$K_1+K_2$$ is an interval.
\end{corollary}
\begin{remark}
The conditions in Corollary \ref{cor2} imply that $\tau(K_1)\tau(K_2)\geq 1$.  It is easy to find some Cantor sets, under the condition $\tau(K_1)\tau(K_2)< 1$, such that $f(K_1,K_2)$ does not contain some interiors, see for instance in Corollary \ref{Main1} and the remarks below.
By the Newhouse's thickness theorem, if $\tau(K_1)\tau(K_2)\geq 1$, then $K_1+K_2$ is an interval. However, under the condition $\tau(K_1)\tau(K_2)\geq 1$, we may not have  that $$K_1\cdot K_2=\{xy:x\in K_1, y\in K_2\}$$ is still an interval.
A  simple example is the middle-third Cantor set, denoted by $C$. The thickness of $C$ is $1$. We have
$$C+C=[0,2].$$ However,   $C\cdot C\subset [0,1/3]\cup [4/9,1]$, which  yields that $C\cdot C$ is not an interval. Therefore, for a general $f$, if we want $f(K_1,K_2)$ to be some interval, we may expect more strong conditions on $f$ besides $\tau(K_1)\tau(K_2)\geq 1$.
\end{remark}
\begin{corollary}\label{cor3}
Let $K_1$ and $K_2$ be two Cantor sets with convex hull $[0,1]$.     If $\tau(K_1)\tau(K_2)>1$, then there  are uncountably many  nonlinear functions
$f(x,y)\in \mathcal{C}^{1}$ such that
$$f(K_1,K_2)$$ is an interval.
\end{corollary}
The condition on partial derivatives in Theorem \ref{Main} can be weakened  when we consider some homogeneous self-similar sets. Indeed, the thickness gives little information about the relation between  gaps and bridges. If we elaborately analyze  their relation, we may obtain more delicate result.
For instance,
with a similar discussion as Theorem \ref{Main}, we may prove the following result.
\begin{corollary}\label{Main1}
Let $K_{\lambda}$ be the attractor of the IFS
$$\{f_1(x)=\lambda x,f_2(x)=\lambda x+1-\lambda, 0<\lambda<1/2 \}.$$
 Suppose that $f(x,y)\in \mathcal{C}^1$ is a continuous function defined on $\mathbb{R}^2$.
 If for any $(x,y)\in [0,1]^2$, we have
$$\dfrac{1-2\lambda}{\lambda}\leq \left|\dfrac{\partial_x f}{\partial_y f}\right|\leq  \dfrac{1}{1-2\lambda},$$ then
$$f(K_{\lambda},K_{\lambda})=\left[\min_{(x,y)\in K_{\lambda}  \times K_{\lambda}}f(x,y),\max_{(x,y)\in K_{\lambda}  \times K_{\lambda}}f(x,y)\right].$$
\end{corollary}
\begin{remark}
The conditions in Corollary \ref{Main1} imply that
$$\dfrac{1-2\lambda}{\lambda}\leq \dfrac{1}{1-2\lambda}, \mbox{ i.e. } 1/4\leq \lambda<1/2.$$ This condition is natural as for any
$f(x,y)\in \mathcal{C}^1$ and  $0< \lambda<1/4$, we have
$$\dim_{H}(f(K_{\lambda}, K_{\lambda}))\leq \dim_{H}(K_{\lambda}\times K_{\lambda})=\dfrac{2\log 2}{-\log\lambda}<1.$$
In other words, if $0< \lambda<1/4$, then $f(K_{\lambda}, K_{\lambda})$ cannot be an interval.

Note that $\tau(K_{\lambda})=\dfrac{\lambda}{1-2\lambda}<1 \mbox{ if }0<\lambda<1/3$. For this case, the Newhouse's thickness theorem does not offer any information for  $f(K_{\lambda},K_{\lambda})$. Moreover, by  \cite[Theorem 2.4]{Astels}, $\gamma(K_{\lambda})=\dfrac{\tau(K_{\lambda})}{\tau(K_{\lambda})+1}=\dfrac{\lambda}{1-\lambda}$, we cannot make use of Astels'  result to consider whether $f(K_{\lambda},K_{\lambda})$ is an interval as for $1/4<\lambda<1/3 $ we have  $2\gamma(K_{\lambda})<1$.
In fact, for the sum of two Cantor sets, the Newhouse's thickness theorem and Astels' thickness theorem are exactly the same.
We mention some related work. In \cite{Pourbarat},  Pourbarat proved under some assumptions that
$$g_1(K_{\lambda_1})+g_2(K_{\lambda_2})=\{g_1(x)+g_2(y):x\in K_{\lambda_1}, y\in K_{\lambda_2}\}$$ contains an interval, where $g_1,g_2\in \mathcal{C}^1$. In Corollary \ref{cor2}, we prove  under some conditions that $f(K_1, K_2)$ is an interval for general Cantor sets.
\end{remark}
In \cite{Takahashi}, Takahashi asked what is the topological structure of $K_{\lambda_1}\cdot K_{\lambda_2}$. He also posed the question for the multiple product of some $K_{\lambda_i}$. In fact, we can simultaneously  consider multiplication and division on $K_{\lambda}$.  We  partially answer his questions as follows.
\begin{corollary}\label{Multiplication}
Let $\{K_{\lambda_i}\}_{i=1}^{d}$  be  self-similar sets with $0<\lambda_i<1/2, i=1,\cdots,d.$ If for any $1\leq i\leq d-1$
\begin{equation*}
\left\lbrace\begin{array}{cc}
\dfrac{1-2\lambda_i}{\lambda_i}\leq 1-\lambda_d\\
\dfrac{1}{1-\lambda_i}\leq \dfrac{\lambda_d}{1-2\lambda_d},
 \end{array}\right.
\end{equation*}
then
$$\Pi_{i=1}^d K_{\lambda_i}^{\epsilon_i}=\{\Pi_{i=1}^{d}x_i^{\epsilon_i}:x_i\in K_{\lambda_i},\epsilon_i\in\{-1,1\}, x_i\neq 0 \mbox{ if } \epsilon_i=-1\}=U,$$
where $$U=\bigcup_{k_1, k_2,\cdots, k_d\in \mathbb{N}}\lambda_1^{\epsilon_1k_1}\lambda_2^{\epsilon_2k_2}\cdots \lambda_d^{\epsilon_dk_d}[\delta, \eta]\cup\{0\},$$
$$\delta= \Pi_{\epsilon_i=1}(1-\lambda_i), \eta=\Pi_{\epsilon_i=-1}(1-\lambda_i)^{-1}.$$
\end{corollary}
\begin{remark}
To avoid triviality, in the definition of $\Pi_{i=1}^d K_{\lambda_i}^{\epsilon_i}$, we assume that
there exist some $1\leq i,j\leq d$ such that  $\epsilon_i=-1, \epsilon_j=1$.  For this case, $\Pi_{i=1}^d K_{\lambda_i}^{\epsilon_i}$ contains $0$. If  $\epsilon_i=-1,  $ for any $1\leq i\leq d$, then $\Pi_{i=1}^d K_{\lambda_i}^{\epsilon_i}$ does not contain  $0.$
\end{remark}
We may find more similar  conditions, as  in the above corollary, which  allow us to describe the structure of  $\Pi_{i=1}^d K_{\lambda_i}^{\epsilon_i}$. Note that
 $$\Pi_{i=1}^d K_{\lambda_i}^{\epsilon_i}=\bigcup_{k_1, k_2,\cdots, k_d\in \mathbb{N}}\lambda_1^{\epsilon_1k_1}\lambda_2^{\epsilon_2k_2}\cdots \lambda_d^{\epsilon_dk_d}(\Pi_{i=1}^d ( \widetilde{K_{\lambda_i}})^{\epsilon_i})\cup\{0\}.$$
 where each $\widetilde{K_{\lambda_i}}$ is the right similar copy of $K_i.$
 The above result only investigates  $f(x_1,\cdots, x_d)$ on  $ \widetilde{K_{\lambda_1}}\times  \cdots \times\widetilde{K_{\lambda_d}}$, see the details in the proof. Indeed, we may decompose each $ \widetilde{K_{\lambda_i}}$ into two sub self-similar sets, and analyze the partial derivatives on these sub similar sets. We leave these considerations to the reader.

 The following result indicates  that the multiplication and division on some self-similar sets may  simultaneously reach their maximal ranges.
\begin{corollary}\label{Division}
Let $K$ be the attractor of the following IFS
$$\{f_1(x)=\lambda_1 x, f_2(x)=\lambda_2x +1-\lambda_2, 0<\lambda_2\leq \lambda_1<1, \lambda_1+\lambda_2<1\}.$$
Then the following conditions are equivalent:
\begin{itemize}
\item [(1)] $$K\cdot  K= \{x  \cdot y:x,y\in K\}=[0,1];$$
\item [(2)] $\lambda_1\geq (1-\lambda_2)^2;$
\item [(3)] $$K\div K=\left\{\dfrac{x}{y}:x,y\in K, y\neq 0\right\}=\mathbb{R}.$$
\end{itemize}
\end{corollary}
Finally, we give an application to the continued fractions with  restricted partial quotients.
We first give some basic definitions. Let $m\in \mathbb{N}_{\geq 2}$. Define
$$F(m)=\{[t,a_1,a_2,\cdots ]: t\in \mathbb{Z}, 1\leq a_i\leq m \mbox{ for }i\geq 1\},$$
where$$ [t,a_1,a_2,\cdots ]=t+\dfrac{1}{a_1+\dfrac{1}{a_2+\dfrac{1}{\cdots}}}. $$
For each $l\in \mathbb{N}^{+}$, define
 \begin{eqnarray*}
G(l)&=& \{[t,a_1,a_2,\cdots ]: t\in \mathbb{Z}, a_i\geq l \mbox{ for }i\geq 1\}\\
&& \cup\{[t,a_1,a_2,\cdots, a_k ]: t,k,\in \mathbb{Z},  k\geq 0, \mbox{ and } a_i\geq l \mbox{ for }1\leq i\leq k\}.
\end{eqnarray*}

Generally, let $B$ be a finite digits set, denote by $F(B)$ the set of points which have an infinite continued fraction expansions with all partial quotients, except possibly the first, members of $B$.  When $B$ is infinite, then we define $F(B)$ in a similar way ($F(B)$ also includes some real numbers with finite continued expansions). For more detailed introduction, see \cite{Astels,KCor}. Let $F_t(B)$ be a subset of  $F(B)$ with the first partial quotient $t\in \mathbb{Z}$. Therefore, $F(B)=\cup_{t\in \mathbb{Z}}F_t(B).$ With this notation, we have
$$F(m)=\cup_{t\in \mathbb{Z}}F_t(B)$$ where $B=\{1,2,\cdots,m\}.$
It is not difficult to calculate the Newhouse thickness of $F_t(B)$, see for instance \cite[Lemma 4.3, Lemma 4.4]{Astels}.

The main motivation why we consider continued fractions with restricted partial quotients  is due to some well-known results.  Hall \cite{Hall} proved that
$$F(4)+F(4)=\{x+y:x,y\in F(4) \}=\mathbb{R}.$$ Divi\v{s} \cite{Divi} showed that Hall's result is sharp in some sense as $$F(3)+F(3)\neq \mathbb{R}.$$ Here we use one simple fact, i.e. $F(n)\subset F(n+1)$ for any $n\geq 2$.  Hlavka  \cite{Hlavka} generalized Hall's result and proved that
$$F(3)+F(4)=\mathbb{R}, F(2)+F(7)=\mathbb{R}, F(2)+F(4)\neq \mathbb{R}.$$
Astels \cite{Astels} showed
$$F(2)\pm F(5)=\mathbb{R}, F(3)- F(4)=\mathbb{R},F(2)- F(4)\neq \mathbb{R},F(3)-F(3)\neq \mathbb{R}.$$
All of the above equations are linear, i.e. the associated function $$f(x,y)=x\pm y$$ is linear.  It is natural to ask can we obtain similar results for some non-linear functions.
Note that  $F_t(B), t\in \mathbb{Z}$ is a Cantor set. Therefore, by Theorem \ref{Main} and the thickness of $F_t(B)$ we can obtain some nonlinear results for the arithmetic on $F(B)$ if we appropriately  choose the functions $f(x,y)$.
We only  give the following equations. The reader may find more similar identities which can represent real numbers.
\begin{corollary}\label{continued fraction}
$$ F^3(7)\pm F(7)=\mathbb{R}, (C+1)^2+2 F(6)=\mathbb{R}, f(K_1, K_2, K_3)=\mathbb{R},$$ where
$$K_1=K_2=C+1, K_3=F(6), f(x,y,z)=0.1x+xy+z,$$
and  $C$ is the middle-third Cantor set.
\end{corollary}
This paper is arranged as follows. In Section 2, we prove the main results of this paper. In Section 3, we give some identities which can represent real numbers. Finally, we give some remarks and questions.
\section{Proofs of main results}
\subsection{Proof of Corollary \ref{cor2}}
We first prove Corollary \ref{cor2}. The proof of  general result, i.e. Theorem \ref{Main}, depends on Corollary \ref{cor2}.

 Clearly,
$$f(K_1,K_2)\subset H.$$ To prove $f(K_1,K_2)=H,$ we suppose on the contrary that $f(K_1,K_2)\neq H$, then we shall find some contradictions. Therefore, we finish the proof of Corollary \ref{cor2}.

\noindent If $$f(K_1,K_2)\neq H,$$ then we can find some  $z\in H$ such that $\Phi_z$ does not intersect with $K_1\times K_2,$ where
$$\Phi_z=\{(x,y)\in [0,1]^2: f(x,y)=z\}.$$
By virtue of  the continuity of $f$, it follows that $\Phi_z$ is a compact set.
Note that   $\Phi_z$ can be covered by  countably many strips  $\Gamma$ of the form
$$\{(x,y):x\in O, y\in [0,1]\}   \mbox{ or }  \{(x,y):y\in O, x\in [0,1]\},$$ where $O$'s are the deleted  open intervals when we construct $K_i, i=1,2$.  By the compactness of $\Phi_z$, we may find finitely many strips from the above covering, i.e. $$\Phi_z\subset \cup_{i=1}^{n}\Gamma_i.$$ By the construction of Cantor sets (we mainly use the fact that the deleted open intervals are pairwise disjoint) and the continuity of  $f(x,y)$, it follows that  for any $1\leq i\leq n-1$,
$\Gamma_i$ is perpendicular to $\Gamma_{i+1}$.

Suppose that $\Gamma_{\min}$ is the strip which has minimal width (every strip has length $1$)  among $\cup_{i=1}^{n}\Gamma_i$. For  every $\Gamma_i, 1\leq i\leq n$, we denote its width  by $L_i.$ Then we have the following  lemma.
	\begin{figure}
	\begin{tikzpicture}[scale=4]
	\draw[->](-0.1,0)--(0.9,0)node[left,below,font=\tiny]{$x$};
	\draw[->](0,-0.1)--(0,0.9)node[right,font=\tiny]{$y$};
\draw[](0,0)--(0,0)node at (-0.05,-0.05){$o$};
\draw[domain=0:0.62] plot(\x,0.22+0.3*\x*\x) node at (0.9,0.35){$f(x,y)=z$};
	\draw[](0,0.3)--(0.8,0.3)node at (-0.10,0.25){$L_{min}$};
	\draw[](0,0.3)--(0.8,0.3)node at (0.56,-0.05){$L_{1}$};
	\draw[](0,0.3)--(0.8,0.3)node at (0.2,0.16){$\rho$};
	\draw[](0,0.2)--(0.8,0.2);
	\draw[](0.45,0)--(0.45,0.8);
	\draw[](0.65,0)--(0.65,0.8);
	\draw[](0,0.3)--(0.8,0.3)node at (-0.05,0.16){$A$};
	\draw[](0,0.3)--(0.8,0.3)node at (0.40,0.16){$B$};
	\draw[](0,0.3)--(0.8,0.3)node at (0.40,0.34){$C$};
	\draw[](0,0.3)--(0.8,0.3)node at (-0.05,0.34){$D$};
	\end{tikzpicture}
\begin{tikzpicture}[scale=4]
	\draw[->](-0.1,0)--(0.9,0)node[left,below,font=\tiny]{$x$};
	\draw[->](0,-0.1)--(0,0.9)node[right,font=\tiny]{$y$};
%\draw[thick] plot [smooth] coordinates{(A) (B) (C) (D)};
% \draw (2.5,4.7) parabola bend (2  .8,2) (6.7,0.75);
\draw (0.8,0.52) .. controls (0.54,0.38) and (0.29,0.6) .. (0.24,0.19) ;
\draw[](0,0)--(0,0)node at (-0.05,-0.05){$o$};
\draw[](0.2,0)--(0.2,0.8);
\draw[](0.35,0)--(0.35,0.8);
\draw[](0.7,0)--(0.7,0.8);
\draw[](0.87,0)--(0.87,0.8);
\draw[](0,0.4)--(0.87,0.4);
\draw[](0,0.5)--(0.87,0.5);
\draw[](0,0.4)--(0.87,0.4)node at (-0.10,0.45){$L_{min}$};
	\draw[](0,0)--(0,0)node at (0.3,-0.05){$L_{2}$};
	\draw[](0,0)--(0,0)node at (0.5,-0.07){$L_{2,3}$};
	\draw[](0,0)--(0,0)node at (0.39,0.37){$A$};
\draw[](0,0)--(0,0)node at (0.66,0.37){$B$};
\draw[](0,0)--(0,0)node at (0.66,0.54){$C$};
\draw[](0,0)--(0,0)node at (0.39,0.54){$D$};
	\draw[](0.2,0)--(0.2,0.8)node at (0.8,-0.05){$L_{3}$};
\draw[](0,0)--(0,0)node at (1.1,0.45){$f(x,y)=z$};
\end{tikzpicture}
\begin{tikzpicture}[scale=3.4]
	\draw[->](-0.1,0)--(1.1,0)node[left,below,font=\tiny]{$x$};
	\draw[->](0,-0.1)--(0,1.1)node[right,font=\tiny]{$y$};
\draw[](0,0)--(0,0)node at (-0.05,-0.05){$o$};
\draw[](0.3,0)--(0.3,1);
\draw[](0.4,0)--(0.4,1);
\draw[](0,0.2)--(1,0.2);
\draw[](0,0.4)--(1,0.4);
\draw[](0,0.7)--(1,0.7);
\draw[](0,0.9)--(1,0.9);
\draw[](0,0)--(0,0)node at (0.35,-0.06){$L_{min}$};
\draw[](0,0)--(0,0)node at (0.6,0.3){$L_{4}$};
\draw[](0,0)--(0,0)node at (0.6,0.8){$L_{5}$};
\draw[](0,0)--(0,0)node at (0.6,0.55){$L_{4,5}$};
\draw (0.5,0.8) .. controls (0.25,0.8) and (0.4,0.1) .. (0.23,0.25) ;
\draw[](0,0)--(0,0)node at (0.65,1){$f(x,y)=z$};
\end{tikzpicture}
\caption{}
	\end{figure}
\begin{lemma}\label{hor}
The strip $\Gamma_{\min}$ does not  parallel with the $x$-axis.
\end{lemma}
\noindent\textbf{Proof:}
We  prove this lemma for three cases. Firstly, if the strip $\Gamma_{\min}$ is parallelling with the $x$-axis, and it is closest to the origin.
Then  by the  implicit function theorem and the minimal width of  $\Gamma_{\min}$ (we denote its width by $L_{min}$),   there exists some $(x_0,y_0)\in \Gamma_{\min}$ such that
$$\left|\dfrac{dy}{dx}|_{(x_0,y_0)}\right|=\left|\dfrac{\partial_x f|_{(x_0,y_0)}}{\partial_y f|_{(x_0,y_0)}}\right|<\dfrac{L_{min}}{\rho}\leq \dfrac{L_{min}}{L_1 \tau(K_1)}\leq \dfrac{1}{ \tau(K_1)},$$ see the first graph of Figure 1. This contradicts to  the condition in Corollary \ref{cor2}. Secondly, if  the strip $\Gamma_{\min}$ is parallelling with the $x$-axis, and it is closest to the line $y=1$, then we may find a similar contradiction as the first case.
Finally, suppose   the strip $\Gamma_{\min}$ is parallelling with the $x$-axis, and there is at least one parallelling strip below and above $\Gamma_{\min}$, respectively. Let $\Gamma_{2}$ and $\Gamma_{3}$ be two strips that are perpendicular to $\Gamma_{\min}$ such that $\Phi_z$ enters and leaves the $\Gamma_{\min}$.  The entrance point is in $\Gamma_{2}$ while the leaving point is in  $\Gamma_{3}$. Let $L_{2,3}$ be the distance between  $\Gamma_{2}$ and $\Gamma_{3}$.
Then by the implicit function theorem   there exists some $(x_0,y_0)\in \Gamma_{\min}$ such that
$$\left|\dfrac{dy}{dx}|_{(x_0,y_0)}\right|=\left|\dfrac{\partial_x f|_{(x_0,y_0)}}{\partial_y f|_{(x_0,y_0)}}\right|<\dfrac{L_{min}}{L_{2,3}}\leq\dfrac{L_{min}}{\min\{L_2, L_3\}\tau{(K_1)}} \leq \dfrac{1}{ \tau(K_1)},$$ see the second graph of Figure 1.
This  contradicts to the assumption of  Corollary \ref{cor2}.  Hence, we have proved Lemma \ref{hor}.

Similarly, we can prove the following lemma.
\begin{lemma}\label{ver}
The strip $\Gamma_{\min}$ does not  parallel with the $y$-axis.
\end{lemma}
\noindent \textbf{Proof:}
Suppose that  $\Gamma_{\min}$ is parallelling with the $y$-axis.
We  prove this lemma in   three cases which are similar to Lemma \ref{hor}. For simplicity, we only prove the following case.

Suppose there is at least one parallelling strip located on the left  and right of $\Gamma_{\min}$, respectively. Then we let $\Gamma_{4}$ and $\Gamma_{5}$ be two strips that are perpendicular to $\Gamma_{\min}$ such that the $\Phi_z$ enters and leaves the $\Gamma_{\min}$. The entrance point is in $\Gamma_{4}$  and  the leaving point is in  $\Gamma_{5}$. Denote by $L_{4,5}$ the distance between $\Gamma_{4}$ and $\Gamma_{5}$, see the third graph of Figure 1.  Therefore, by the implicit function theorem again, there exists some $(x_0,y_0)\in\Gamma_{\min} $ such that
$$\left|\dfrac{dy}{dx}|_{(x_0,y_0)}\right|=\left|\dfrac{\partial_x f|_{(x_0,y_0)}}{\partial_y f|_{(x_0,y_0)}}\right|>\dfrac{L_{4,5}}{L_{\min}}\geq \dfrac{\min\{L_{4}, L_{5}\}\tau(K_2)}{L_{\min}}\geq \tau(K_2) .$$
This is a contradiction.

\textbf{Proof of  Corollary \ref{cor2}}
 Corollary \ref{cor2} follows from Lemmas \ref{hor}
and \ref{ver}.
\subsection{Proof of Theorem \ref{Main}}
Now, we prove Theorem \ref{Main}. The main idea is exactly the same as Corollary \ref{cor2}.
Firstly, we clearly have $$f(K_1,\cdots,  K_d)\subset H. $$ If
$$f(K_1,\cdots,  K_d)\subsetneq H, $$ then  there exists some $w\in H$ such that
the hypersurface $f(x_1,\cdots,  z)=w$ does not intersect with $$K_1\times K_2\times\cdots \times K_d.$$
Now  we construct the following set
$$\Psi_w=\{(x_1,x_2,\cdots, z)\in [0,1]^d: f(x_1,\cdots,  z)=w\}.$$ It is a compact set by the continuity of $f$.
Hence, we can find   finitely many $d$-dimensional cubes of the form $$\Lambda=\Delta_1\times \Delta_2\times \cdots \times  \Delta_d\subset [0,1]^d,$$
 such that there is a unique $\Delta_i \subsetneq [0,1]$ is an open interval for some $1\leq i\leq d$, and the rest $\Delta_j=[0,1], j\neq i$. We call each $\Delta_i, 1\leq i\leq d$ an edge of $\Delta_1\times \Delta_2\times \cdots \times  \Delta_d$.   For simplicity, we call the edge  which is not equal to $[0, 1]$ the axis edge.
Without loss of generality, we may assume that  $\Psi_w\subset \cup_{i=1}^{n}\Lambda_i$,  $\Lambda_i $ is  perpendicular to $  \Lambda_{i+1} $ for $1\leq i\leq n-1$, i.e. the  axis edges of $\Lambda_i $ and $\Lambda_{i+1} $ have different subscripts.  Let $\Lambda_{min}$ be the cube with minimal length,  i.e.  one  edge of $\Lambda_{min}$ has minimal length among $\cup_{i=1}^{n}\Lambda_i$.  We shall prove that the above covering, i.e. $\cup_{i=1}^{n}\Lambda_i$, does not exists. Therefore, we prove the desired result.

\noindent Let $$\Lambda_{min}=[0,1]^{i-1} \times (p_i, q_i) \times[0,1]^{d-i}, (p_i, q_i) \subsetneq [0,1], 1\leq i\leq d.$$
Suppose $1\leq i\leq d-1$, for the hypersurface $$f(x_1,\cdots,  z)=w,$$ we fix $x_j, j\neq i,d $ (we let $x_d=z$).
Therefore the hypersurface $$f(x_1,\cdots,  z)=w$$
can be covered by  $\Omega_i\cup \Omega_d$, where
$$\Omega_i=\cup\{x_1\}\times\{x_2\}\times  \cdots \times \{x_{i-1}\}\times (p_i,q_i)\times  \{x_{i+1}\}\times  \cdots \times  \{x_{d-1}\}\times [0,1],$$
and $$\Omega_d=\cup\{x_1\}\times\{x_2\}\times  \cdots \times \{x_{i-1}\}\times [0,1]\times  \{x_{i+1}\}\times  \cdots \times  \{x_{d-1}\}\times (p_d,q_d).$$
Here the unions in the above equations mean finite (by the compactness of $\Psi_w$) deleted open intervals when we construct $K_i$ and $K_d.$

Since we fix $x_j, j\neq i,d $, it follows that the hypersurface $f(x_1,\cdots,  z)=w$  becomes a curve on a plane which is a translation of the $x_iOz$ plane.  We let this curve be $\Upsilon$.  By the above discussion, we have  $\Upsilon\subset \Omega_i\cup \Omega_d.$
Nevertheless, by Lemmas \ref{hor} and  \ref{ver}, $\Upsilon$ cannot be in $\Omega_i\cup \Omega_d$, which is a contradiction.
If $i=d$, then $$\Lambda_{min}=[0,1]^{d-1} \times (p_d, q_d).$$ For this case, we can also prove similarly as above, and obtain a contradiction.
 Hence, we finish the proof.
\subsection{Proof of  Corollary \ref{cor3}}
Note that
$$\tau(K_1)\tau(K_2)>1\Leftrightarrow \dfrac{1}{\tau(K_1)}<\tau(K_2).$$
If $\tau(K_2)>\tau(K_1)>1$, then there exist some $\alpha, \beta \in \mathbb{R}^{+} $ such that
$$\tau(K_2)> 1+2\alpha, \tau(K_1)>1+2\beta. $$
Now, we let $f(x,y)=\alpha x^2+\beta y^2+x+y$. Since the convex hull of $K_i, i=1,2$ is $[0,1]$, it follows that the conditions in Corollary \ref{cor2} are satisfied.  Therefore, $f(K_1, K_2)$ is an interval.

If $\tau(K_2)>1>\tau(K_1)$, then we can find some $\gamma, \zeta\in \mathbb{R}^{+}$ such that
$$\dfrac{1}{\tau(K_1)}<\dfrac{\gamma}{2+\zeta},\dfrac{\gamma+2}{\zeta}<\tau(K_2) .$$ We let $$f(x,y)=x^2+y^2+\gamma x+\zeta y.$$
It is easy to check the conditions in Corollary \ref{cor2}. Hence,  $f(K_1, K_2)$ is an interval.
\subsection{Proof of  Corollary \ref{Main1}}
 The proof is almost the same as the proof of Corollary \ref{cor2}. We only need to prove Lemma \ref{ver} under the assumption
$$ \left|\dfrac{\partial_x f}{\partial_y f}\right|\leq \dfrac{1}{1-2\lambda}.$$ For simplicity, we only prove the first case of Lemma \ref{ver}. We still use the  terminology  of Lemma \ref{ver}, and the  third graph of Figure 1.
 By Lemma \ref{hor} $\Gamma_{min}$ cannot parallel with $x$-axis. By the minimality of $\Gamma_{min}$, we have
 $$\lambda\min\{L_{4}, L_{5}\}\geq  L_{min}.$$
 Therefore, by the implicit function theorem,  there exists some $(x_0,y_0)\in\Gamma_{\min} $ such that
$$\left|\dfrac{dy}{dx}|_{(x_0,y_0)}\right|=\left|\dfrac{\partial_x f|_{(x_0,y_0)}}{\partial_y f|_{(x_0,y_0)}}\right|>\dfrac{L_{4,5}}{L_{\min}}\geq \dfrac{\min\{L_{4}, L_{5}\}\tau(K_{\lambda})}{L_{\min}} \geq \dfrac{\tau(K_\lambda)}{\lambda}=\dfrac{1}{1-2\lambda} .$$
This is a contradiction.
\subsection{Proof of  Corollary \ref{Multiplication}}
Let $f(x_1,\cdots ,x_d)=\Pi_{i=1}^{d}x_i^{\epsilon_i}$. It is easy to check that
$$\left|\dfrac{\partial_{x_i}f}{\partial_{x_d}f}\right|=\left|\dfrac{x_d}{x_i}\right|, 1\leq i\leq d-1.$$
Note that
 $$\Pi_{i=1}^d K_{\lambda_i}^{\epsilon_i}=\bigcup_{k_1, k_2,\cdots, k_d\in \mathbb{N}}\lambda_1^{\epsilon_1k_1}\lambda_2^{\epsilon_2k_2}\cdots \lambda_d^{\epsilon_dk_d}(\Pi_{i=1}^d ( \widetilde{K_{\lambda_i}})^{\epsilon_i})\cup\{0\},$$
 where $\widetilde{K_{\lambda_i}}$ is the right similar copy of $K_i.$ Note that the convex hull of $$  \widetilde{K_{\lambda_1}}\times  \widetilde{K_{\lambda_2}}\times \cdots  \times\widetilde{K_{\lambda_d}}$$ is
 $$V=[1-\lambda_1,1]\times [1-\lambda_2,1]\times \cdots\times [1-\lambda_d,1]. $$
Therefore,   for any $(x_1,\cdots, x_d)\in V$, we have
$$1-\lambda_d\leq \left|\dfrac{x_d}{x_i}\right|\leq \dfrac{1}{1-\lambda_i}.$$ Then by the following conditions
\begin{equation*}
\left\lbrace\begin{array}{cc}
\dfrac{1-2\lambda_i}{\lambda_i}\leq 1-\lambda_d\\
\dfrac{1}{1-\lambda_i}\leq \dfrac{\lambda_d}{1-2\lambda_d},
 \end{array}\right.
\end{equation*}
we clearly have
$$\tau(K_{\lambda_i})^{-1}\leq \left|\dfrac{x_d}{x_i}\right|\leq \tau(K_{\lambda_d}), 1\leq i\leq d-1.$$
Now, Corollary \ref{Multiplication} follows from Theorem \ref{Main}.
\subsection{Proof of Corollary \ref{Division}}
 We first prove $(1)\Rightarrow (2)$. This is clear as
 $$K\subset [0,\lambda_1]\cup [1-\lambda_2,1]\Rightarrow K\cdot K\subset [0,\lambda_1]\cup [(1-\lambda_2)^2,1].$$
    Now, we prove that $(2)\Rightarrow (1)$. Let $f(x,y)=xy$. First, we have the following equation:
    $$K\cdot K=\cup_{i=0}^{\infty}\lambda_1^{i}(f_2(K)\cdot f_2(K))\cup\{0\}.$$
    The convex hull of $f_2(K)$ is $[1-\lambda_2,1]$.  Hence, we consider the partial derivatives of $f$ on $[1-\lambda_2,1]^2$. It is easy to calculate that
    $$1-\lambda_2\leq \left|\dfrac{\partial_x f}{\partial_y f}\right|=\left|\dfrac{y}{x}\right|\leq \dfrac{1}{1-\lambda_2} \mbox{ for any } (x,y)\in [1-\lambda_2,1]^2.$$
    Note that $\lambda_1\geq (1-\lambda_2)^2$ is equivalent to $\dfrac{1}{1-\lambda_2}\leq  \tau(K)=\dfrac{\lambda_2}{1-\lambda_1-\lambda_2}$. Therefore, by Corollary \ref{cor2}, $$f_2(K)\cdot f_2(K)=[(1-\lambda_2)^2,1].$$
    Since $\lambda_1\geq (1-\lambda_2)^2$, it follows that
    $$K\cdot K=\cup_{i=0}^{\infty}\lambda_1^{i}(f_2(K)\cdot f_2(K))\cup\{0\}=[0,1].$$
    Now we prove $(3)\Rightarrow (1)$. Note that
    $$K\div K= \bigcup_{n=-\infty}^{+\infty}\lambda_1^n\dfrac{f_2(K)}{f_2(K)}\cup \{0\}\subset\bigcup_{n=-\infty}^{+\infty}\lambda_1^n\left(\left[1-\lambda_2, \dfrac{1}{1-\lambda_2}\right]\right)\cup \{0\}. $$
    Therefore, if $\lambda_1<(1-\lambda_2)^2$, then
    $$ \left[1-\lambda_2, \dfrac{1}{1-\lambda_2}\right]\bigcap \left[\lambda_1(1-\lambda_2), \dfrac{\lambda_1}{1-\lambda_2}\right]=\emptyset.$$ In other words, $$K\div K\neq \mathbb{R}.$$
    Finally, we prove $(1)\Rightarrow (3)$. This step is  almost the same  as $(2)\Rightarrow (1)$ in terms of  the equation
  $$K\div K= \bigcup_{n=-\infty}^{+\infty}\lambda_1^n\dfrac{f_2(K)}{f_2(K)}\cup \{0\}.$$
 \section{Some identities}
In this section, we mainly prove Corollary \ref{continued fraction}.
It is easy to calculate $$\tau(F_t(7))=(42 + 24\sqrt{77})/91, t\in \mathbb{Z},$$  see \cite[Lemma 4.3, Lemma 4.4]{Astels}.
  Therefore, by Corollary \ref{cor2}, it follows that $$F_1^3(7)+F_1(7)=\left[(\dfrac{7+\sqrt{77}}{14})^3+\dfrac{7+\sqrt{77}}{14}, (\dfrac{-5+\sqrt{77}}{2})^3+\dfrac{-5+\sqrt{77}}{2}\right].$$ Moreover, it is easy to check that
$$(F_1^3(7)+F_i(7))\cap (F_1^3(7)+F_{i+1}(7))\neq \emptyset , i\in \mathbb{Z}.$$
Therefore, $$F^3(7)+F(7)=\mathbb{R}.$$ Similarly, we can prove $$F^3(7)-F(7)=\mathbb{R}.$$
For the second identity, we first note that
$$\dfrac{1}{2}(C+1)^2+ F(6)=\mathbb{R}\Leftrightarrow (C+1)^2+2 F(6)=\mathbb{R},$$
where $$(C+1)^2+2 F(6)=\{x^2+2y:x\in C+1, y\in F(6)\}.$$
Hence, we only need to prove $$\dfrac{1}{2}(C+1)^2+ F(6)=\mathbb{R}.$$
 Let
$$f(x,y)=\dfrac{1}{2}x^2+y, x\in C+1, y\in F_0(6).$$ Then by Corollary \ref{cor2}, we have
$$\dfrac{1}{2}(C+1)^2+ F_0(6) $$ is an interval and
$$\left(\dfrac{1}{2}(C+1)^2+ F_i(6)\right)\cap \left(\dfrac{1}{2}(C+1)^2+ F_{i+1}(6)\right)\neq \emptyset , i\in \mathbb{Z}.$$ As such,
$$\dfrac{1}{2}(C+1)^2+ F(6)=\mathbb{R}.$$
Finally, we consider the function
$$f(x,y,z)=0.1x+xy+z, x,y\in C+1, z\in F_0(6).$$
By Theorem \ref{Main}, we have that
$$f(C+1,C+1,F_0(6))=\left[1.1+\dfrac{-3+\sqrt{15}}{6}, 1.2+\sqrt{15}\right].$$ Moreover,
$$f(C+1,C+1,F_i(6))\cap f(C+1,C+1,F_{i+1}(6))\neq \emptyset, i\in \mathbb{Z}.$$
Therefore, we have
$$f(C+1,C+1,F(6))=\mathbb{R}.$$
\section{Final remarks and some problems}
Although in Theorem \ref{Main}, we give a sufficient condition under which the continuous image of $f $ is a closed interval, there are many problems left. We list some problems as follows.
\begin{itemize}
\item [(1)]For a given $E\subset \mathbb{R}^d$, define a continuous function $g:\mathbb{R}^d\to \mathbb{R}^d$. It would be interesting to consider when $g(E)$ contains an interior or $g(E)$ is exactly some convex hull.
\item [(2)] In Theorem \ref{Main}, we do not know whether for two concrete sets, the lower and upper bounds  of the ratio of partial derivatives can be improved.
\item [(3)] In Theorem \ref{Main}, we only consider the first order partial derivatives. Can we give a similar nonlinear version of Theorem \ref{Main} using higher orders of  partial derivatives.
\item[(4)] In Corollary \ref{Division}, we find an example such that the resonant maximum for the multiplication and division occurs.  It would be interesting to find more sets which have this resonant phenomenon. Moreover, we may consider the resonant phenomenon for other arithmetic operation such as sum of squares and sum of cubes. These questions are motivated by the representations of real numbers from number theory.
\item[(5)] Given two Cantor sets $K_1$ and $K_2$ with $\tau(K_1)\tau(K_2)<1$, can we find some sufficient conditions such that $f(K_1, K_2)$ is still an interval.
\item[(6)] Given two Cantor sets $K_1$ and $K_2$, we do not know when $f(K_1, K_2)$ is a union of finitely many closed intervals.
    \item[(7)] The Newhouse's thickness, in some sense, is rough. As it gives a rough relation between gaps and bridges. It is deserved  to define   a finer thickness. Under the new thickness, we may partially improve Theorem \ref{Main}.
\end{itemize}

%\bibliographystyle{plain}
%\bibliography{OnUnivoquePointsForSelfSimilarSets}

\begin{thebibliography}{99}

\bibitem{Astels}  Steve Astels. \newblock Cantor sets and numbers with
restricted partial quotients. \newblock {\em Trans. Amer. Math. Soc.},
352(1):133--170, 2000.

\bibitem{Tyson}
Jayadev~S. Athreya, Bruce Reznick, and Jeremy~T. Tyson.
\newblock Cantor set arithmetic.
\newblock {\em Amer. Math. Monthly}, 126(1):4--17, 2019.

\bibitem{BJJ}
Taras Banakh, Eliza Jab{\l}o\'{n}ska, and Wojciech Jab{\l}o\'{n}ski.
\newblock The continuity of additive and convex functions which are upper
  bounded on non-flat continua in {$\Bbb R^n$}.
\newblock {\em Topol. Methods Nonlinear Anal.}, 54(1):247--256, 2019.

\bibitem{KCor} Karma Dajani and Cor Kraaikamp.
\newblock {\em Ergodic theory
of numbers}, volume~29 of \emph{Carus Mathematical Monographs}. \newblock %
Mathematical Association of America, Washington, DC, 2002.

\bibitem{Divi}
Bohuslav Divi\v{s}.
\newblock On the sums of continued fractions.
\newblock {\em Acta Arith.}, 22:157--173, 1973.

\bibitem{Fraser}
Jonathan M. Fraser, Douglas C. Howroyd and  Han Yu.
\newblock Dimension growth for iterated sumsets.
\newblock {\em Math. Z.}, 293(3-4):1015--1042, 2019.

\bibitem{Feng2020}
Dejun Feng and Yufeng Wu. \newblock    On arithmetic sums of fractal sets in $\mathbb{R}^d$.  \newblock {\em
 arXiv:2006.12058}, 2020.


\bibitem{Gu}Jiangwen Gu, Kan Jiang, Lifeng Xi, and Bing Zhao. \newblock Multiplication on uniform $\lambda$-Cantor set. \newblock {\em
 arXiv:1910.08303}, 2019.


 \bibitem{Hall}
Marshall Hall, Jr.
\newblock On the sum and product of continued fractions.
\newblock {\em Ann. of Math. (2)}, 48:966--993, 1947.


\bibitem{Hlavka}
James~L. Hlavka.
\newblock Results on sums of continued fractions.
\newblock {\em Trans. Amer. Math. Soc.}, 211:123--134, 1975.



\bibitem{Yoccoz}
Carlos Gustavo~T. de~A.~Moreira and Jean-Christophe Yoccoz.
\newblock Stable intersections of regular {C}antor sets with large {H}ausdorff
  dimensions.
\newblock {\em Ann. of Math. (2)}, 154(1):45--96, 2001.






\bibitem{Newhouse}
Sheldon~E. Newhouse.
\newblock The abundance of wild hyperbolic sets and nonsmooth stable sets for
  diffeomorphisms.
\newblock {\em Inst. Hautes \'{E}tudes Sci. Publ. Math.}, (50):101--151, 1979.

\bibitem{Palis}
Jacob Palis and Floris Takens.
\newblock {\em Hyperbolicity and sensitive chaotic dynamics at homoclinic
  bifurcations}, volume~35 of {\em Cambridge Studies in Advanced Mathematics}.
\newblock Cambridge University Press, Cambridge, 1993.
\newblock Fractal dimensions and infinitely many attractors.

\bibitem{Pourbarat}
Mehdi.Pourbarat.
\newblock On the arithmetic difference of middle cantor sets.
\newblock {\em Discrete and Continuous Dynamical Systems.}, 38(9):4259-4278, 2018.


\bibitem {HS}Hugo Steinhaus. \newblock Mowa W{\l }asno\'{s}\'{c} Mnogo\'{s}ci
Cantora. \newblock {\em Wector, 1-3. English translation in: STENIHAUS,
H.D.} 1985.



\bibitem{Takahashi}
Yuki Takahashi.
\newblock Products of two Cantor sets.
\newblock {\em Nonlinearity}, 30(5):2114--2137, 2017.

\bibitem{Yu}Han Yu. \newblock Fractal projections with an application in number theory. \newblock {\em
 arXiv:2004.05924}, 2020.
\end{thebibliography}

  \section*{Acknowledgements}
  This work is
supported by K.C. Wong Magna Fund in Ningbo University.
This work is also supported by National Natural Science Foundation of China with No.  11701302, and by Zhejiang Provincial Natural Science Foundation of China with
No.LY20A010009.

\end{document}